\documentclass[10pt,a4paper]{article}
\usepackage{amsfonts}
\usepackage {graphicx}
\usepackage{amsmath}
\usepackage{amssymb}
\usepackage{amsthm}
\usepackage{natbib}
\numberwithin{equation}{section}
\numberwithin{figure}{section}

\newtheorem{theorem}{\sc Theorem}[section]
 
\newtheorem{proposition}[theorem]{\sc Proposition} 
\newtheorem{lemma}[theorem]{\sc Lemma} 
 
\newtheorem{definition}[theorem]{\sc Definition}

\title{The Exact Value for  European Options on a Stock Paying a Discrete
Dividend}

\author{Jo\~{a}o Amaro de Matos\dag,  Rui Dil\~{a}o\ddag  \ and  Bruno Ferreira\ddag}

\begin{document}
\date{}

\maketitle

\begin{center}
\dag Faculdade de Economia, Universidade Nova de Lisboa,
Rua Marqu\^{e}s de Fronteira, 20, 1099-038 Lisbon, Portugal
\end{center}

\begin{center}
\ddag Nonlinear Dynamics Group, Instituto Superior T\'{e}cnico,
Av. Rovisco Pais, 1049-001 Lisbon, Portugal
\end{center}

\begin{center}
amatos@fe.unl.pt; rui@sd.ist.utl.pt; bruno.ferreira@mckinsey.com
\end{center}

\begin{abstract}
In the context of a Black-Scholes
economy and with a no-arbitrage argument, we derive  arbitrarily accurate lower and upper bounds for the value of European
options  on  a stock paying a discrete dividend. 
Setting the option price error below the smallest monetary unity, both bounds coincide, and we
obtain the exact value of the option.
\end{abstract}

\vskip 2 true cm



\maketitle

\section{Introduction}

In the seminal paper of Black and Scholes (1973), the problem of valuing a
European option was solved in closed form. Among other things, their
result assumes that the stochastic process associated to
the underlying asset is a geometric Brownian motion, not allowing for
the payment of discrete dividends. Yet the majority of stocks on which
options trade do pay dividends.

Merton (1973) was the first to relax the no-dividend assumption, allowing
for a deterministic dividend yield. In this case, he showed that  European
options can be  priced in the context of a Black-Scholes economy, with
either a continuous dividend yield or a discrete dividend proportional to
the stock price. However, when the dividend process is discrete and does not
depend on the stock level, the simplicity of the Black-Scholes model breaks
down.

Let $S_{t}$ denote the value of the underlying asset at time $t$, and  let $T$ be  the maturity time of the option. When the
risky asset pays a dividend $D$ at time $\tau <T$, a jump of size $D$ in the
value process happens at that point in time. 
The stock price process is discontinuous at $t=\tau $
and is no more a geometric Brownian motion in the time interval $[0,T]$.

The
standard approximation procedure for valuing European options written on
such a risky asset, first informally suggested by Black (1975), considers
a  Black-Scholes formula, where the initial price of the underlying
stock $S_{0}$ is replaced by its actual value less the present value ($PV$) of the
dividends  ($Div$), 
\[
S_{0}\to S_{0}^{\ast }=S_{0}-PV(\hbox{Div}) 
\]
This adjustment is made to evaluate the option at any point in time before
$\tau $.  After the payment of dividends, there is no need for further adjustments.
In this approximation, the input in the Black-Scholes formula  is the value of
the (continuous) stochastic process, 
\[
S_{t}^{\ast }=\left\{ 
\begin{array}{l}
S_{t}-De^{-r\left( \tau -t\right) },\quad t<\tau \\ 
S_{t},\quad t\geq \tau
\end{array}
\right.
\]
where $r$ is the risk-free rate. 

For $t<\tau $, the discontinuous stock price process $S_{t}$ can thus be seen as the sum of
two components ($S_{t}=S_{t}^{\ast }+De^{-r\left( \tau
-t\right) }$). One riskless component, $De^{-r\left( \tau -t\right) }$,
corresponding to the known dividends during the life of the option, and a
continuous risky component $S_{t}^{\ast }$. At any
given time before $\tau $, the riskless component   is the present value of the dividend discounted
at the present at the risk-free rate. For any time after $\tau $ until the
time the option matures, the dividend will have been paid and the riskless
component will no longer exist. We thus have $S_{T}=S_{T}^{\ast }$ and, as
pointed out by Roll (1977), the usual Black-Scholes formula is correct to
evaluate the option  only   if  $S_{t}^{\ast }$ 
follows a geometric Brownian motion. In that case, we would use in the Black-Scholes
formula $S_{0}^{\ast }$ for the initial value, together with the volatility
of the process $S_{t}^{\ast }$, followed by the risky component of
the underlying asset.

If we assume that $S_{t}^{\ast }$\ follows a geometric Brownian motion,
a simple application of It\^{o} Lemma shows that the original stock price
process $S_{t}$ does not follow a geometric Brownian motion in the time interval $\left[
0,\tau \right[ $. On the other hand, under the Black-Scholes assumption that
 $S_{t}$ follows a geometric Brownian motion in $\left[ 0,\tau \right[ $,
the risky component $S_{t}^{\ast }$ follows a continuous
process that is not a geometric Brownian motion in $\left[0,\tau \right[ $. 
Therefore, the standard procedure described  above must be
seen  as an approximation to the true value of such calls under the
Black-Scholes assumption. As argued by Bos and
Vandermark (2002),  this assumption is  typically underlying the
intuition of traders, but the approximation is sometimes bad. In
fact, as noticed in the early papers about option pricing 
(Cox and Ross, 1976; Merton, 1976a; Merton, 1976b), the correct specification
of the stochastic process followed by the value of the underlying stock is
of prime importance in option valuation.

The deficiency of this standard procedure is reported in Beneder and Vorst
(2001). Using Monte Carlo simulation methods, these authors calculate the values of call options 
under the Black-Scholes assumption, and compare them with
the values obtained with the approach just described. Reported errors are up
to  $9.4\%$. They also find that the standard procedure above usually
undervalues the options. For these reasons, Beneder and Vorst (2001) propose a
different approximation,  trying to improve the standard procedure by 
adjusting the volatility of the underlying asset. This approach consists in
modifying the variance of the returns by  a weighted average of an
adjusted and an unadjusted variance, where the weighting depends on the
time $\tau $ of the dividend payment. Performing much better than the
former approximation, this method still does not allow the  control of the  errors
committed for the given parameters of the economy. Analogously, Frishling
(2002) warns on the mispricing risk due to the use of an incorrect  underlying 
stochastic process. This discussion is followed by a series of
recent papers suggesting different approximations that better match
numerical results (Bos and Vandermark, 2002;  Bos \textit{et al}, 2003). 
More recently, 
Haug  \textit{et al} (2003) discuss this problem.
However, as these authors claim, ``[i]n the case of European
options, the above techniques are \textit{ad hoc}, but the job gets done (in
most cases) when the corrections are properly carried out''. 

The development of these approximations enhance two important aspects.
First, they are not exact, and it is not possible to control the
error with respect to the correct value of the option.
Second, there are numerical procedures to estimate the value of these options, as for example,
Monte-Carlo simulation methods. However, this method is time consuming and provides a 
convergence of statistical nature.

The purpose of this paper is to derive a closed form for the exact value of
European options on a stock paying a discrete dividend, in the context of a Black-Scholes
economy. We obtain an exact
result and we need not to rely on \textit{ad hoc} assumptions. 

This paper is organized as follows. In Section 2, an integral representation
for the value of  European options written on an asset paying a discrete
dividend is obtained, and the convexity properties of the solutions of the Black-Scholes equation are
derived. In section 3, we construct functional upper and lower bounds for
the integral representation of the value of an option. These bounds follow from a
convexity property of the solutions of the Black-Scholes equation. Theorem \ref{main} is the main result of this paper and gives the  algorithmic procedure to determine the price of European options  on a stock
paying a discrete dividend.
In section 4,  numerical examples are analyzed and we discuss the advantages 
of the proposed method.
In section 5, we summarize the main conclusions of the paper.

\section{Valuation of European options on a stock paying a discrete dividend}

In this section, following a standard procedure to derive the Black-Scholes formula 
(Wilmott, 2000),  we derive an integral representation for the
value of a European option written on an asset paying a known discrete
dividend. 

We consider a European call option with maturity time $T$ and
strike price $K$. This call option is written on an underlying asset with
value $S_{t}$, with stochastic
differential equation, 
\[
dS_{t}=\mu S_{t}dt+\sigma S_{t}dW_{t}
\]
where $\mu $ and $\sigma $ are the drift and volatility of the underlying
asset. The quantity $W_{t}$ is a continuous and normally distributed
stochastic process with mean zero and variance $t$. Under these conditions, 
the underlying asset with value $S_{t}$  follows a geometric Brownian motion.
We also assume a
risk-free asset with constant rate of return $r.$

In the context of the Black-Scholes economy, the value $V$ of an option is
dependent of the time $t$ and of the price of the underlying asset $S$. Under
the absence of arbitrage opportunities (Wilmott, 2000; 
Bj\"{o}rk, 1998), it follows that $V(S,t)$ obeys the Black-Scholes equation, 
\begin{equation}\label{eq:pde}
{\frac{{\partial V}}{{\partial t}}}+{\frac{1}{2}}\sigma ^{2}S^{2}{\frac{{ 
\partial ^{2}V}}{{\partial S^{2}}}}+rS{\frac{{\partial V}}{{\partial S}}} -rV=0   
\end{equation}

The Black-Scholes equation is a quasi-linear parabolic partial differential
equation, with $S\geq 0$, and $t\geq 0$. To determine the solutions of the
Black-Scholes equation, we introduce the new variables, 
\begin{equation}
\label{eq:newvariables}
\left\{
\begin{array}{l}
\theta =T-t   \nonumber \\
x =\log S+\left( r-\frac{\sigma ^{2}}{2}\right) \left( T-t\right) \nonumber \\
\end{array}
\right.
\end{equation}
together with the new function $\varphi (x,\theta )=e^{r(T-t)}V(S,t)$. In
the new coordinates (\ref{eq:newvariables}), the Black-Scholes equation (\ref
{eq:pde}) becomes the diffusion equation, 
\begin{equation}
{\frac{{\partial \varphi }}{{\partial \theta }}}={\frac{1}{2}}\sigma ^{2}{
\frac{{\partial ^{2}\varphi }}{{\partial x^{2}}}}  \label{eq:diffusion}
\end{equation}
where $x\in \mathbf{R}$ and $\theta \geq 0$. If $\theta =0$, by (\ref
{eq:newvariables}), we have $\varphi (x,0)=V(S,T)$, and $\varphi
(x,T)=e^{rT}V(S,0)$. Therefore, by (\ref{eq:newvariables}), the forward solution in the time $\theta $
of the diffusion equation relates  with the
backward solution in the time $t$ of the Black-Scholes equation (\ref{eq:pde}%
). The Black-Scholes problem for the price of a call option is to determine
the option value at time $t=0$ whose value at maturity time $T$ is, 
\begin{equation}\label{eq:fronteira}
V\left( S,T\right) =\max \{0,S-K\}  
\end{equation}
Therefore, due to the change of coordinates (\ref{eq:newvariables}), the
call option solution of  the Black-Scholes equation (\ref{eq:pde}) is equivalent to an initial
value problem for the diffusion equation.

Suppose now an initial data problem for the diffusion equation (\ref
{eq:diffusion}), $\varphi (x,\theta =0)=f(x)$. Under these
conditions, the general solution of (\ref{eq:diffusion}) is (Folland, 1995), 
\begin{equation}
\varphi \left( {x,\theta }\right) ={\frac{1}{{\sigma \sqrt{2\pi \theta }}}}
\int_{-\infty }^{\infty }f(y){\exp }\left[ {{-{\frac{{\left( {x-y}\right)
^{2}}}{{2\sigma ^{2}\theta }}}}}\right] {dy}  \label{eq:soldiff}
\end{equation}
and the solution of the Black-Scholes equation for a call option is, 
\begin{equation}
V\left( {S}{,0}\right) =e^{-rT}\varphi \left( {x,T}\right) ={\frac{
e^{-rT}}{{\sigma \sqrt{2\pi T}}}}\int_{-\infty }^{\infty }V(e^{y},T){\exp }
\left[ {{-{\frac{{\left( {x-y}\right) ^{2}}}{{2\sigma ^{2}T}}}}}\right] {dy}
\label{eq:solBS}
\end{equation}
This integral can be easily calculated to obtain the usual Black-Scholes formula
(Black and Scholes, 1973; Wilmott, 2000).

For a dividend distribution at some time $\tau \in (0,T)$, the
Black-Scholes formula is no longer true, since, during the life time of
the option, the value of the underlying
asset does not follow a geometric Brownian motion. However, if we take the time intervals, $I_{1}=[0,\tau \lbrack 
$ and $I_{2}=[\tau ,T]$, the value of the underlying asset follows a
geometric Brownian motion in each interval $I_{1}$ and $I_{2}$, and, at time $t=\tau $, it has a
jump equal to the dividend $D$.

Before considering this  case, we proceed with some properties
of the  solutions (\ref{eq:soldiff}) and (\ref{eq:solBS}) of the
diffusion and of the Black-Scholes equations.

\begin{definition}
A real valued function $f\left( x \right)$, with $x \in \mathbf{R}$, is convex if,
for every $x_1 , x_2 \in \mathbf{R}$, 
\[
f\left( {{\frac{{x_1 + x_2 } }{2}}} \right) \le {\frac{1 }{2}}\left( {
f\left( {x_1 } \right) + f\left( {x_2 } \right)} \right)
\]
\end{definition}

A simple property of convex functions is that, if the real-valued functions $f$ and $g$ are both convex,
and $g$ is increasing, then $f(g(x))$ is also convex.

\begin{proposition}
\label{convex} Let $f(x)$ the initial data function of a well-posed
diffusion equation problem, and suppose that $f(x)$ is non-negative and convex. 
Then, for fixed $\theta $, the solution $\varphi (x,\theta )$ of the diffusion equation 
 is also convex. Moreover, if  $f(x)$ is an increasing function, then, for fixed $\theta $,
 $\varphi (x,\theta )$ is also increasing.
\end{proposition}

\emph{Proof.} Suppose that the solution (\ref{eq:soldiff}) of the diffusion equation  (\ref{eq:diffusion}) is well defined
(Folland, 1995).
By (\ref{eq:soldiff}), with $z=y-x$, we have, 
\[
\varphi \left( {x,\theta }\right) ={\frac{1}{{\sigma \sqrt{2\pi \theta }}}} 
\int_{-\infty }^{\infty }f(z+x){\exp }\left( {{-{\frac{{z^{2}}}{{2\sigma
^{2}\theta }}}}}\right) {dz}
\]
As, by hypothesis, $f(x)$ is convex, then, for every $z\in \mathbf{R}$, 
\[
f\left[ {{\frac{{\left( {x_{1}+z}\right) +\left( {x_{2}+z}\right) }}{2}}}
\right] =f\left( z+{\frac{x_{1}+x_{2}}{2}}\right) \leq {\frac{1}{2}}\left[ { 
f\left( {z+x_{1}}\right) +f\left( {z+x_{2}}\right) }\right]
\]
and, as $f(x)$ is non-negative, 
\begin{eqnarray*}
\varphi \left( {\frac{x_{1}+x_{2}}{2},s}\right) &=&\frac{1}{{\sigma \sqrt{ 
2\pi \theta }}}\int_{-\infty }^{\infty }f\left( z+{\frac{x_{1}+x_{2}}{2}} 
\right) {\exp }\left( {{-{\frac{{z^{2}}}{{2\sigma ^{2}\theta }}}}}\right) {dz 
} \\
&\leq &{\frac{1}{2}}\left[ \varphi \left( x_{1},\theta \right) +\varphi
\left( x_{2},\theta \right) \right]
\end{eqnarray*}
and so $\varphi (x,\theta )$ is also convex. 
Assuming now that $f(x)$ is increasing, we have that $f(x_2)\ge f(x_1)$, whenever $x_2>x_1$. Then,
for every $z \in \mathbf{R}$, we have, 
$f(z+x_2)\ge (z+x_1)$, and, by (\ref{eq:soldiff}), the last assertion of the proposition follows. \hfill $\Box$

\bigskip

As (\ref{eq:fronteira}) is a convex function in $S$, Proposition \ref{convex}
implies that the backward solution (\ref{eq:solBS}) of the Black-Scholes equation (\ref
{eq:pde}) is also a convex function. 

Suppose now that a dividend on the underlying asset is distributed at time 
$t=\tau $. We denote this dividend by $D$. According to the classical
solution of the Black-Scholes equation (Wilmott, 2000), the price of the option just after
the distribution of dividends at time $t=\tau $ is, 
\begin{equation}
\label{eq:value+}
V\left( {S_{+},\tau }\right) =S_{+}\,N\left( {d+\sigma \sqrt{T-\tau }}%
\right) -Ke^{-r\left( {T-\tau }\right) }N\left( {d}\right) 
\end{equation}
where, 
\[
d={\frac{{\ln S_{+}-\ln K+\left( r-{\frac{1}{2}}\sigma ^{2}\right) \left(
T-\tau \right) }}{{\sigma \sqrt{T-\tau }}} }
\]
and $S_{+}$ denotes the value of the underlying asset just after the
dividend distribution. The function $N(\cdot )$ is the cumulative distribution
function for the normal distribution with mean zero and unit variance. By
Proposition \ref{convex}, the function $V(S_{+},\tau )$ is convex. Note that, the solution
(\ref{eq:value+}) is given by, $V\left( {S_{+},\tau }\right) =e^{-r\left( {T-\tau }\right) } \phi(x, T-\tau)$,
and is directly calculated from (\ref{eq:solBS}) and (\ref{eq:fronteira}).

The approach taken here to value an option is equivalent (see, among others,
Cox and Ross, 1976; Harrison and Krebs, 1979) to write this value at any
point in time as the expected discounted payoff of the option at maturity 
$T, $ under the so-called risk-neutral probability measure. Hence, knowing
beforehand the amount to be distributed as dividend, the value of the option
is not supposed to jump at $\tau $. In other words, the payment of known
dividends $D$ at a known point in time $\tau $ does not affect the
expectations at time $\tau $ about the final payoff of the option at
maturity $T$, and the value of the option is continuous at $\tau $\footnote{According to Wilmott, 2000, the jump condition on the asset price is known  {\it a priori}, implying that there is no surprise
in the fall of the stock price. Therefore, in order to avoid arbitrage opportunities, the value of the option should not change  across the dividend date. This is a no-arbitrage argument.} (Wilmott, 2000, pp. 129-131).
Going backward in time, the value of the underlying asset jumps  from 
$S_{+}$ to $S_{-}=S_{+}+D$, where $S_{-}$ is the value of the underlying asset just
before the dividend distribution. As 
$V(S_{+},\tau)=V(S_{-},\tau )$, by (\ref{eq:value+}), the price of the option just before
the distribution of dividends at time $t=\tau $ is, 
\begin{equation}
\label{eq:value-}
 V(S_{-},\tau )=\left\{
\begin{array}{ll}
(S_{-}-D)N( {\bar{d}+\sigma \sqrt{T-\tau }}) -Ke^{-r(T-\tau)}N( {\bar{d}})  &\hbox{if}\  S_{-}>D \\ 
0 &\hbox{if}\   S_{-}\leq D \\ 
\end{array}
\right.
\end{equation}
where, 
\begin{equation}
{\bar{d}}={\frac{{{\ln (S_{-}-D)-\ln {K}+\left( {r-{\frac{1}{2}}\sigma ^{2}} 
\right) \left( {T-\tau }\right) }}}{{\sigma \sqrt{T-\tau }}}}   \label{eq:d-}
\end{equation}

\begin{figure} [htbp]
\centerline{\includegraphics[width=8 true cm]{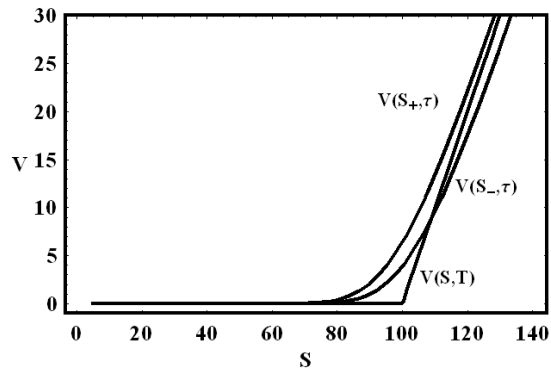}}
\caption{Option values $V({S_+,\tau})$, $V({S_-,\tau})$ and $V(S,T)$ as a function of  the value $S$
of the underlying asset. Parameter values are: 
$\mu =0.01$, $\sigma =0.2$, $r=0.03$, $K=100$, $D=5$, $T=1$ and $\tau=0.5$.}
\label{fig:BS}
\end{figure}

In Fig.  \ref{fig:BS}, we plot $V({S_{+},\tau })$, $V({S_{-},\tau })$ and $V(S,T)$ as a function of $S$.
The  functions $V({S_{+},\tau })$, $V({S_{-},\tau })$ and $V(S,T)$ are convex.

To calculate the value of a call option as a function of the actual price ($t=0$) of
the underlying asset, we must introduce the change of coordinates (\ref
{eq:newvariables}) into (\ref{eq:value-}) and integrate as in (\ref
{eq:solBS}). By (\ref{eq:solBS}) and (\ref{eq:value-}), it follows that the
time-zero value of a European option written on an asset paying dividend $D$
at time $t=\tau $ is given by, 
\begin{equation}
 V\left( {S,0}\right) =e^{-r\tau }\varphi \left( {x,\tau }\right) ={\frac{%
e^{-r\tau }}{{\sigma \sqrt{2\pi \tau }}}}\int_{-\infty }^{\infty }V\left[ {%
S_{-}(y),\tau }\right] {\exp }\left[ {{-{\frac{{\left( {x-y}\right) ^{2}}}{{%
2\sigma ^{2}\tau }}}}}\right] {dy}  \label{eq:solBSwithD}
\end{equation}
which has no simple representation in terms of tabulated functions. By  
Proposition \ref{convex}, $V(S,0)$ is also convex.

\section{Accurate bounds for $V( S,0)$}

As it is difficult to determine a close form for the integral representation of
the option's value (\ref{eq:solBSwithD}) in terms of tabulated functions,
to estimate the value  $V\left( {S,0}\right) $, we use
the convexity property of $V({S_{-},\tau })$ and its asymptotic behavior as $%
S_{-}\rightarrow \infty $.

\begin{lemma}
\label{asymptotic}If $K>0$, then, in the limit $S_{-}\rightarrow \infty $, $%
V({S_{-},\tau })$ is asymptotic to the line $V=(S_{-}-D)-Ke^{-r(T-\tau )}$,
and $V({S_{-},\tau })\geq (S_{-}-D)-Ke^{-r(T-\tau )}$.
\end{lemma}

\emph{Proof.}
\noindent In the limit $S_{-}\rightarrow \infty $, ${\bar{d}}\rightarrow
\infty $, and $N({\bar{d}})\rightarrow 1$. Hence, by (\ref{eq:value-}), 
$V({S_{-},\tau })$ is asymptotic to the line $V_1=(S_{-}-D)-Ke^{-r(T-\tau )}$. To
prove the second part of the lemma, first note that, if  $V_1=(S_{-}-D)-Ke^{-r(T-\tau )}\le 0$, then $S_{-}\leq
D+Ke^{-r(T-\tau )}$. As $V({S_{-},\tau })$ is non-negative, if $S_{-}\leq
D+Ke^{-r(T-\tau )}$, then $V({S_{-},\tau })\ge V_1$. 

Suppose now that 
$S_{-}>D+Ke^{-r(T-\tau )}$.   
By hypothesis, we  assume  that there exists some $S_{-}={\bar{S}}
$ such that, $V({{\bar{S}},\tau })=({\bar{S}}-D)-Ke^{-r(T-\tau )}$, and $V({{\bar{S}},\tau })>0$. 
By (\ref{eq:value-}) and (\ref{eq:d-}), we then have, 
\[
Ke^{-r(T-\tau )}={\frac{N\left[ {\bar{d}}({\bar{S}})+\sigma \sqrt{T-\tau }%
\right] -1}{N\left[ {\bar{d}}({\bar{S}})\right] -1}}({\bar{S}}-D)
\]
As $(S_{-}-D)>  Ke^{-r(T-\tau )}$, from the equality above, we obtain,
\[
{\frac{N\left[ {\bar{d}}({\bar{S}})+\sigma \sqrt{T-\tau } 
\right] -1}{N\left[ {\bar{d}}({\bar{S}})\right] -1}}({\bar{S}}-D)=Ke^{-r(T-\tau )}<(S_{-}-D)
\]
 Hence,
\[
N\left[ {\bar{d}}({\bar{S}})+\sigma \sqrt{T-\tau }\right] <N\left[ {\bar{d}}({\bar{S}})\right]
\]
which contradicts the fact that $N(\cdot )$ is a monotonically increasing
function of the argument. Therefore,  the function $V({S_{-},\tau })$ and the line 
$V_1=(S_{-}-D)-Ke^{-r(T-\tau )}$ do not intersect for finite $\bar{S}$. 
As $V({S_{-},\tau })$ is a continuous function of $S_{-}$,
then $V({S_{-},\tau })\ge V_1$ in all the range of $S_{-}$, and the lemma is proved. \hfill
$\Box$

\medskip

To estimate the solution (\ref{eq:solBSwithD}) of the Black-Scholes
equation, we use Proposition \ref{convex} and Lemma \ref{asymptotic} to
construct  integrable upper and lower bound functions of  $V({%
S_{-},\tau })$. This constructions proceeds as follows.

Let us choose a fixed number $S_{-}=S^{\ast }>D$, and divide the interval $%
[D,S^{\ast }]$ into $M\geq 1$ smaller subintervals. The length of the
subintervals is $\Delta S=(S^{\ast }-D)/M$, and their extreme points are
denoted by, 
\[
S_{i}=D+i\,\Delta S,\qquad i=0,\ldots ,M
\]
As the function $V({S_{-},\tau })$ is convex, in each subinterval, the
function $V({S_{-},\tau })$ is bounded from above by the chord that connects
the points $(S_{i},V({S_{i},\tau }))$ and $(S_{i+1},V({S_{i+1},\tau }))$. We
define the constants, 
\[
\alpha _{i}={\frac{M}{S^{\ast }-D}}\left[ V({S_{i},\tau })-V({S_{i-1},\tau })%
\right] ,\qquad i=1,\ldots ,M
\]
where by (\ref{eq:value-}), $V({S_{0},\tau })=0$. Therefore, in each
interval $[S_{i-1},S_{i}]$, the function $V({S_{-},\tau })$ is bounded from
above by the function $f_{i}(S_{-})=\alpha _{i}(S_{-}-S_{i-1})+V({%
S_{i-1},\tau })$.

Let us define the characteristic function of a set $I$ as, $\chi _{I}(x)=1$, if $%
x\in I$, and $\chi _{I}(x)=0$, otherwise. Then, the function $V({S_{-},\tau }%
)$ in the interval $[D,S^{\ast }]$ is approached from above by the piecewise
linear function, 
\begin{equation}
V_{1}^{+}({S_{-},\tau })=\sum_{i=1}^{M}\left[ \alpha _{i}(S_{-}-S_{i-1})+V({%
S_{i-1},\tau })\right] \chi _{\lbrack S_{i-1},S_{i}]}(S_{-})  \label{eq:upb1}
\end{equation}

To extend the bound of $V({S_{-},\tau })$ to $S_{-}>S^{\ast }$, we introduce
the function, 
\begin{equation}
V_{2}^{+}(S_{-},\tau )=\left[ (S_{-}-S^{\ast })+V(S^{\ast },\tau )\right]
\chi _{\lbrack S^{\ast },\infty )}(S_{-})  \label{eq:upb2}
\end{equation}

By Proposition \ref{convex} and Lemma \ref{asymptotic}, for $S_{-}\geq
S^{\ast }$, $V_{2}^{+}(S_{-},\tau )$ is the chord connecting the point $%
(S^{\ast },V(S^{\ast },\tau ))$ to the point at infinity. Therefore, we have
proved the following:

\begin{lemma}
\label{upper} The function $V(S_{-},\tau )$ has the upper bound,  
\[
V(S_{-},\tau )\leq V_{1}^{+}(S_{-},\tau )+V_{2}^{+}(S_{-},\tau ),\ \ \hbox{if}\ \ S_{-}> D
\]
where $V_{1}^{+}$ and $V_{2}^{+}$ are given by (\ref{eq:upb1}) and (\ref
{eq:upb2}), respectively, and the function $(V_{1}^{+}+V_{2}^{+})$ is piecewise linear and 
non-negative. If $S_{-}\le D$, $V(S_{-},\tau )=0$.
\end{lemma}

The construction of a lower bound for (\ref{eq:value-}) follows the same line of reasoning.

In each subinterval $[S_{i-1},S_{i}]\subset \lbrack D,S^{\ast }]$, we can
construct a linear function that bounds from below the function $%
V(S_{-},\tau )$. Due to the convexity of $V(S_{-},\tau )$, we construct the
lower bound through the derivative of $V(S_{-},\tau )$ at the middle point
of each interval $[S_{i-1},S_{i}]$. We then have, 
\begin{equation}
 V_{1}^{-}\left( {S_{-},\tau }\right) =\sum\limits_{i=1}^{M}\left[ {{%
V^{\prime }\left( {S_{i+{\frac{1}{2}}},\tau }\right) \left( {S_{-}-S_{i+{%
\frac{1}{2}}}}\right) +V\left( {S_{i+{\frac{1}{2}}},\tau }\right) }}\right]
\chi _{\lbrack S_{i-1},S_{i}]}(S_{-})  \label{eq:V1-}
\end{equation}
where, 
\[
V^{\prime }\left( {S_{-},\tau }\right) ={\frac{{e^{-{\frac{1}{2}}\left( {%
\overline{d}+\sigma \sqrt{T-\tau }}\right) ^{2}}}}{{\sigma \sqrt{2\pi }\sqrt{%
T-\tau }}}}-{\frac{{K\,e^{-r\left( {T-\tau }\right) }\,e^{-{\frac{1}{2}%
\overline{d}}^{2}}}}{{\sigma \sqrt{2\pi }\sqrt{T-\tau }\left( {S_{-}-D}%
\right) }}}+N\left( {\overline{d}+\sigma \sqrt{T-\tau }}\right)
\]
and ${\overline{d}}$ is given by (\ref{eq:d-}).

To extend the lower bound of $V({S_{-},\tau })$ to $S_{-}>S^{\ast }$, we use
Lemma \ref{asymptotic} to introduce the function, 
\begin{equation}
V_{2}^{-}({S_{-},\tau })=\left[ (S_{-}-D)-Ke^{-r(T-\tau )}\right] \chi
_{\lbrack S^{\ast },\infty )}(S_{-})  \label{eq:V2-}
\end{equation}
By Lemma \ref{asymptotic}, $V_{2}^{-}({S_{-},\tau })$
bounds from below $V(S_{-},\tau )$. Therefore, we have:

\begin{lemma}
\label{lower} The function $V(S_{-},\tau )$ has the lower bound,
\[
V(S_{-},\tau )\geq V_{1}^{-}(S_{-},\tau )+V_{2}^{-}(S_{-},\tau ),\ \ \hbox{if}\ \ S_{-}> D
\]
where $V_{1}^{-}$ and $V_{2}^{-}$ are given by (\ref{eq:V1-}) and (\ref
{eq:V2-}), respectively, and the function $(V_{1}^{-}+V_{2}^{-})$ is piecewise linear and 
non-negative. If $S_{-}\le D$, $V(S_{-},\tau )=0$.
\end{lemma}

Finally, we can state our main result:

\begin{theorem}
\label{main}We consider the Black-Scholes equation (\ref{eq:pde}) together with the terminal condition
(\ref{eq:fronteira}). We assume that $K>0$ and a dividend $D>0$  is payed
at the time $\tau $ with $0<\tau <T$. Let $S=S^{\ast }>D$ be a fixed
constant and let $M\geq 1$ be an integer. Then, the solution of the
Black-Scholes equation with terminal condition (\ref{eq:fronteira}) has the following upper 
and lower bounds: 
\begin{eqnarray*}
V(S,0)\leq V_{S^{\ast },M}^{+}(S,0)&=&\sum\limits_{i=1}^{M}\left\{
\alpha _{i}A_{i}S+e^{-r\tau }\left[ V(S_{i-1},\tau )-\alpha _{i}S_{i-1}%
\right] B_{i}\right\}  \\
&+&SN(d^{\ast })+e^{-r\tau }\left[ V(S^{\ast },\tau )-S^{\ast }\right]
N(d^{\ast }-{\sigma \sqrt{\tau }})
\end{eqnarray*}
and 
\begin{eqnarray*}
V(S,0)\geq V_{S^{\ast },M}^{-}(S,0)&=&S\sum\limits_{i=1}^{M}V^{\prime
}\left( S_{i+\frac{1}{2}},\tau \right) A_{i} \\
&+&e^{-r\,\tau }\sum\limits_{i=1}^{M}\left[ {V}\left( {{{S_{i+{\frac{1}{2}}
},\tau }}}\right) {{-V^{\prime }}}\left( {{{S_{i+{\frac{1}{2}}},\tau 
}}}\right) {{{S_{i+{\frac{1}{2}}}}}}\right] {B_{i}} \\
&+&SN(d^{\ast })-e^{-r\tau }\left( {D+Ke^{-r\left( {T-\tau }\right) }}
\right) N\left( {d^{\ast }-\sigma \sqrt{\tau }}\right)
\end{eqnarray*}
where, 
\begin{eqnarray*}
S_{i} &=&D+{\frac{{S^{\ast }-D}}{M}}i \\
d_{i} &=&{\frac{{\log S-\log S_{i}+(r+{\frac{1}{2}}\sigma ^{2})\tau }}{{%
\sigma \sqrt{\tau }}}} \\
d &=&{\frac{{\log (S-D)-\log K+(r+{\frac{1}{2}}\sigma ^{2})(T-\tau )}}{{%
\sigma \sqrt{T-\tau }}}} \\
d^{\ast } &=&{\frac{{\log S-\log S^{\ast }+(r+{\frac{1}{2}}\sigma ^{2})\tau }%
}{{\sigma \sqrt{\tau }}}} \\
V(S,\tau ) &=&(S-D)N(d)-Ke^{-r(T-\tau )}N(d-{\sigma \sqrt{T-\tau }}) \\
V^{\prime }\left( {S,\tau }\right) &=&N(d)+{\frac{{e^{-{\frac{1}{2}}%
d^{2}}}}{{\sigma \sqrt{2\pi (T-\tau )}}}}-{\frac{{K\,e^{-r\left( {T-\tau }%
\right) }\,e^{-{\frac{1}{2}}\left( {d-\sigma \sqrt{T-\tau }}\right) ^{2}}}}{{%
\sigma \sqrt{2\pi (T-\tau )}\left( {S-D}\right) }}} \\
\alpha _{i} &=&{\frac{M}{{S^{\ast }-D}}}\left[ {V(S_{i},\tau }{%
)-V(S_{i-1},\tau }{)}\right] \\
A_{i} &=&N(d_{i-1})-N(d_{i}) \\
B_{i} &=&N(d_{i-1}-{\sigma \sqrt{\tau }})-N(d_{i}-{\sigma \sqrt{\tau }})
\end{eqnarray*}
and $N(\cdot )$ is the cumulative distribution function for the normal
distribution with mean zero and unit variance.
\end{theorem}

\emph{Proof.}
By Lemmata (\ref{upper}) and (\ref{lower}), 
\[
V_{1}^{-}({S_{-},\tau })+V_{2}^{-}({S_{-},\tau })\leq V({S_{-},\tau })\leq
V_{1}^{+}({S_{-},\tau })+V_{2}^{+}({S_{-},\tau }),\ \ \hbox{if}\ \ S_{-}> D
\]
Multiplying this inequality by the factors as in the integral (\ref{eq:solBSwithD}), and
integrating, we obtain
the estimates of the theorem. \hfill $\Box$

Note that, for $S^{\ast }>D$  fixed,  $\lim_{M\to\infty} V_{S^{\ast },M}^{-}(S,0)\not=
lim_{M\to\infty} V_{S^{\ast },M}^{+}(S,0)$. However, if  $S^{\ast }$ is large enough, both limits 
can be made arbitrarily close. Technically, this is due to the way the exponential term 
in (\ref{eq:solBS}) contributes to the integral.

\section{Calculating the price of a call option on a stock paying a discrete dividend}

Theorem \ref{main} is the necessary tool to determine  the price of a call option when the underlying asset pays a
discrete known dividend before maturity time $T$. In fact, Theorem \ref{main} asserts
that we can always find  upper and a lower bound functions for 
$V(S,0)$, and the bounding functions approach  each other as we increase $M$ and $S^{\ast }$.

\begin{figure} [htbp]
\centerline{\includegraphics[width=14 true cm]{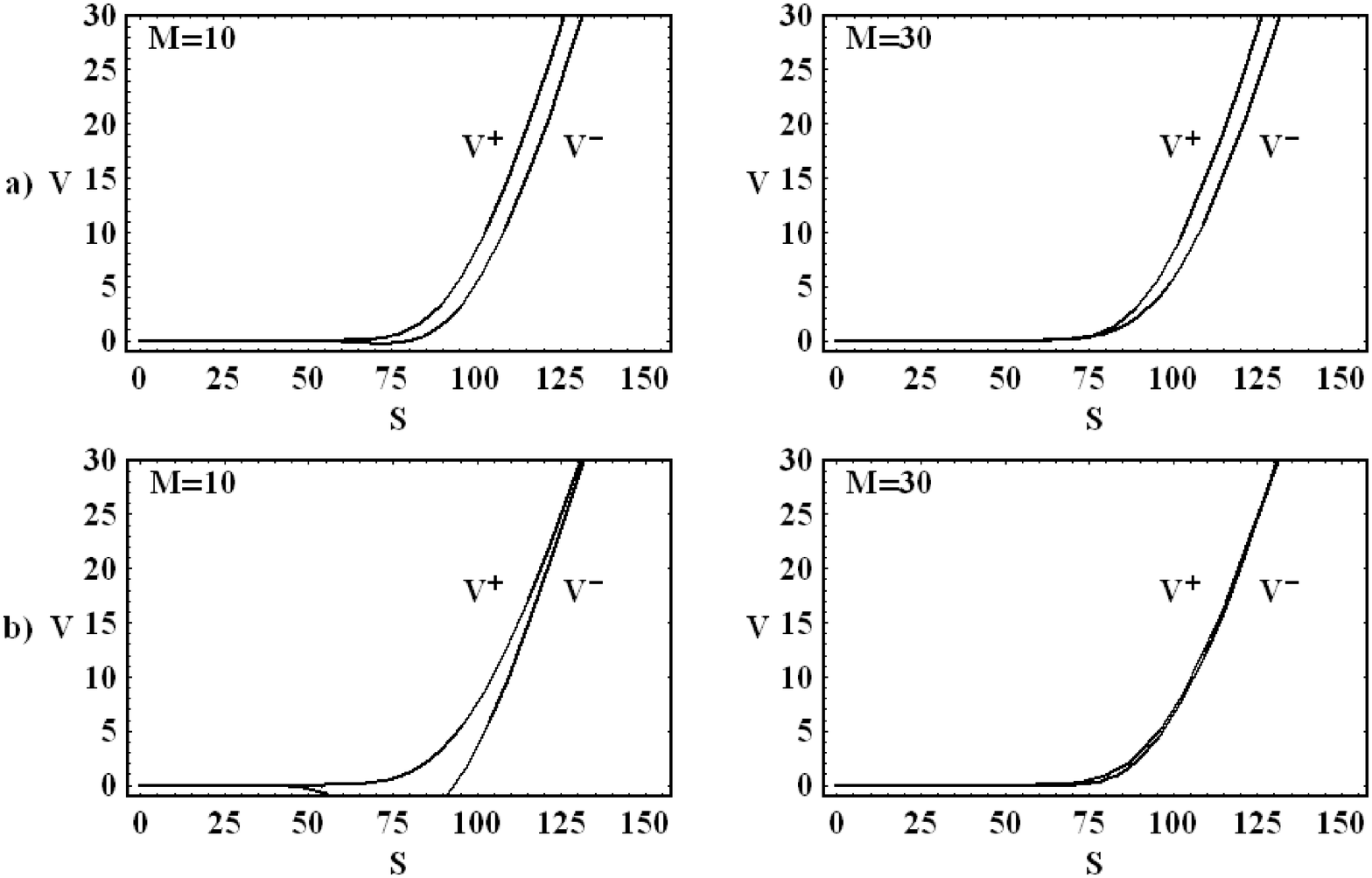}}
\caption{Bounds $V^+_{S^*,M}(S,0)$ and $V^-_{S^*,M}(S,0)$ for $V(S,0)$, calculated
from Theorem \ref{main}, for several values of $S^*$ and $M$. In a) we have chosen  $S^*=D+
Ke^{-r(T-\tau)}=103.5$. In b), $S^*=2(D+ Ke^{-r(T-\tau)})=207.0$. Parameter values are: $\mu =0.01$, 
$\sigma =0.2$, $r=0.03$, $K=100$, $D=5$, $T=1$ and $\tau=0.5$.}
\label{fig:bounds}
\end{figure}

To determine the price of the option, we first choose fixed values for  the 
approximation parameters
$S^{\ast }$ and $M$. If $V_{S^{\ast },M}^{+}(S,0)$ and $V_{S^{\ast
},M}^{-}(S,0)$ differ too much within some fixed precision, we then
increase $S^{\ast }$ and $M$.

\bigskip

\begin{table} 
\caption{\label{tb:t1} Bounds $V_{S^{\ast },M}^{+}(S,0)$ and $V_{S^{\ast },M}^{-}(S,0)$
for $V(S,0)$, calculated from Theorem \ref{main}, for several values of $S^{\ast }$ 
and $M$, and $S=110$. 
The exact
value  $V(S,0)$ has been obtained by the numerical integration of (\ref
{eq:solBSwithD}). The interval error $\varepsilon$ is given by 
(\ref{eq:solBSwithD}). Parameter values are the same as in 
Fig. \ref{fig:bounds}, and we have chosen $S^{\ast }=D+Ke^{-r(T-\protect\tau
)}=103.5$,  $S^{\ast }=1.5(D+Ke^{-r(T-\protect\tau )})=155.3$ and $S^{\ast }=2(D+Ke^{-r(T-\protect\tau )})=207.0$.}

\begin{center}
\item[]\begin{tabular}{ccccccc} 
\hline
S & $S^{\ast }$ & $M$ & $V_{S^{\ast },M}^{-}(S,0)$ & $V(S,0)$ & $V_{S^{\ast
},M}^{+}(S,0)$ & $\varepsilon$\\ 
\hline
110 & 103.5 & 10 & 11.24 & 12.87 & 15.41 & 4.166 \\ 
110 & 103.5 & 50 & 11.61 & 12.87 & 15.35 & 3.739\\ 
110 & 103.5 & 400 & 11.63 & 12.87 & 15.35 & 3.721\\ 
&  &  &  &  &  &\\ 
110 & 155.3 & 10 & 11.39 & 12.87 & 13.20 & 1.807 \\ 
110 & 155.3 & 50 & 12.79 & 12.87 & 12.88 & 0.096\\ 
110 & 155.3 & 200 & 12.87 & 12.87 & 12.87 & 0.006\\ 
110 & 155.3 & 400 & 12.87 & 12.87 & 12.87 & 0.002\\ 
&  &  &  &  &  &\\ 
110 & 207.0 & 10 & 10.64 & 12.87 & 13.45 & 2.813\\ 
110 & 207.0 & 50 & 12.72 & 12.87 & 12.89 & 0.170\\ 
110 & 207.0 & 200 & 12.86 & 12.87 & 12.87 & 0.011\\
110 & 207.0 & 400 & 12.87 & 12.87 & 12.87 & 0.003\\
\hline 
\end{tabular} 
\end{center}
\end{table}

To analyze the convergence of the functional bounds $V^{+}$ and $V^{-}$ to the true price
of a call option,  
we take, as an example, the  parameters: $\mu =0.01$ (drift), 
$\sigma =0.2$ (volatility), $r=0.03$ (interest rate), $K=100$ (strike price), 
$D=5$ (dividend), $T=1$ (expiration time) and $\tau =0.5$ (time of dividend
paying). In Fig. \ref{fig:bounds}, we show $V_{S^{\ast },M}^{+}(S,0)$ and 
$V_{S^{\ast },M}^{-}(S,0)$, for several values of $S^{\ast }$ and $M$, and calculated from Theorem
\ref{main}. Increasing
$M$ and $S^*$, the upper and lower bounds 
$V_{S^{\ast },M}^{+}(S,0)$ and $V_{S^{\ast },M}^{-}(S,0)$ approach each other,
increasing the accuracy to which the functionals bounds approach the option price. To quantify this approximation to
the value of the option, we define the interval error as,
\begin{equation}
\varepsilon=|V_{S^{\ast },M}^{+}(S,0)-V_{S^{\ast },M}^{-}(S,0)|  \label{eq:error}
\end{equation}

In Table \ref{tb:t1}, we compare the values of the upper and lower bounds $V^+_{S^*,M}(S,0)$ and 
$V^-_{S^*,M}(S,0)$, calculated from Theorem \ref{main}, with the exact
value of $V(S,0)$, obtained by the numerical integration of (\ref
{eq:solBSwithD}). We show also the interval error $\varepsilon$ associated to both bounds.
Assuming an interval error below the smallest unit of the monetary currency, for example,
$\varepsilon<10^{-2}$, we obtain the true value of the option. Therefore, for a choice of $S^*$ and $M$
such that $\varepsilon<10^{-2}$, the difference between $V^+_{S^*,M}(S,0)$ and 
$V^-_{S^*,M}(S,0)$, is below the smallest unit of the monetary currency, and the rounded values
of $V^+_{S^*,M}(S,0)$ and  $V^-_{S^*,M}(S,0)$ coincide. This rounded value is the option value
within the chosen monetary accuracy

\begin{figure} [htbp]
\centerline{\includegraphics[width=10 true cm]{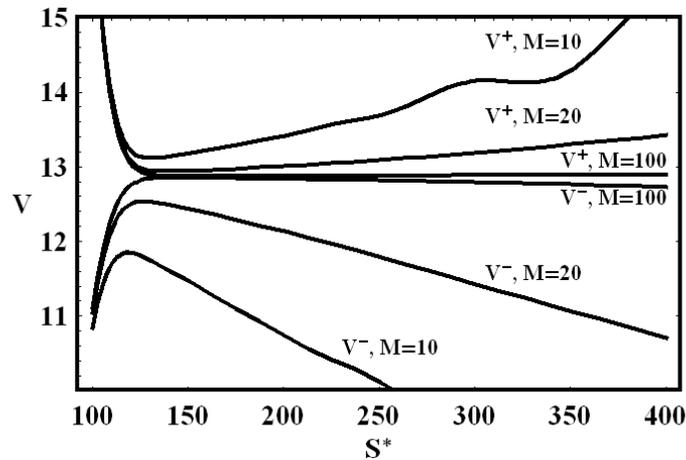}}
\caption{Bounds $V_{S^{\ast },M}^{+}(S,0)$ and $V_{S^{\ast },M}^{-}(S,0)$ as a
function of $S^{\ast }$, for $S=110$ and several values of $M$. The
parameter values are the same as in Fig. \ref{fig:bounds} and Table \ref{tb:t1}. 
}
\label{fig:Md}
\end{figure}

To analyze the global convergence behavior of $V^+_{S^*,M}(S,0)$ and  $V^-_{S^*,M}(S,0)$,
we chose a fixed value of $S$, and we change the approximation 
parameters $S^{\ast }$ and $M$.
In Fig. \ref{fig:Md}, we show $V_{S^{\ast
},M}^{+}(S,0)$ and $V_{S^{\ast },M}^{-}(S,0)$ as a function of $S^{\ast }$,
for several values of $M$. Increasing $M$, the upper and lower bounds of $V(S,0)$ become 
close in a region of the $S^{\ast }$ axis. A choice of $S^{\ast }$ in this region, gives 
better bounds to the value of the option, for lower values of $M$ (Table \ref{tb:t1} and Fig. \ref{fig:Md}). 

For all the examples we have analyzed, 
a good compromise to determine the value of the call option is to choose 
$S^{\ast}=2(D+Ke^{-r(T-\tau )})$. Then, increasing $M$,  
the interval error decreases. Due to the 
fast computational convergence of the expressions in Theorem \ref{main},
bounds with interval error below the smallest unit of the monetary currency
are straightforwardly obtained.

\section{Concluding remarks}

We have obtained an upper and a lower bound for the exact value of a call
option on a stock paying a known discrete dividend at a known future time.
We have assumed the context of a Black-Scholes economy, where, away from the dividend time
paying, the underlying asset price follows a  geometric Brownian motion type stochastic process.
The upper and lower bounds both approach  the exact value of the option
when two parameters are varied. In practical terms, one of these parameters ($S^{\ast }$)
can be fixed to the value, 
$S^{\ast }=2\left( D+Ke^{-r(T-\tau )}\right) $,
where $K$ is the strike, $D$ is the dividend, $\tau $ is the time of paying the
discrete dividend, and $T$ is the length of the contract. Increasing the
second parameter $M$, we obtain bounds for the option value with increasing accuracy. If this accuracy
is below the smallest unit of the  monetary currency, both bounds coincide, and we
obtain the exact value of the option.

The technique used to construct these bounds relies on the convexity
properties of the option value at maturity, and on a property of the
Black-Scholes and diffusion equations that preserves the convexity of
propagated initial conditions. Under this framework, a similar
methodology  can  be used to determine the value of a put
option on a stock paying a known discrete dividend at a known future time.

From the numerical point of view, the technique developed here reduces to
the sum of a few Black-Scholes type terms, whereas numerical Monte Carlo
methods rely on the poor convergence properties determined by the classical
central limit theorem. In our numerical tests for the determination of the exact price of a call
option, the computing time of our technique (using the Mathematica
programming language) is several orders of magnitude faster than the
computing time of finite diferences integration algorithms and of Monte
Carlo methods.

\section*{Acknowledgments}
 We would like to thank Faisal Al-Sharji for the carefull test of the results presented here.
This work has been partially supported by Funda\c c\~ao para a Ci\^encia e a Tecnologia (Portugal), 
under a plurianual funding grant.

\end{document}